\documentclass[sn-mathphys,Numbered]{sn-jnl}

\usepackage{graphicx}%
\usepackage{multirow}%
\usepackage{amsmath,amssymb,amsfonts}%
\usepackage{amsthm}%
\usepackage{mathrsfs}%
\usepackage[title]{appendix}%
\usepackage{xcolor}%
\usepackage{textcomp}%
\usepackage{manyfoot}%
\usepackage{booktabs}%
\usepackage{algorithm}%
\usepackage{algorithmicx}%
\usepackage{algpseudocode}%
\usepackage{listings}%
\usepackage{amsfonts}
\usepackage{xfrac}
\usepackage{bm}
\usepackage{graphicx}
\usepackage{epstopdf}
\usepackage{todonotes}
\usepackage{booktabs}
\usepackage[T1]{fontenc}
\usepackage{tikz}
\usepackage[caption=false,font=footnotesize]{subfig}
\usepackage{pgfplots}
\usepackage{tikz-network}
\pgfplotsset{compat=newest}
\pgfplotsset{plot coordinates/math parser=false}
\usepackage{url}
\usepackage{color}
\usepackage{adjustbox}
\usepackage{enumitem}
\usepackage{placeins} 

\usepackage{mathtools}
\mathtoolsset{showonlyrefs}

\theoremstyle{thmstyleone}%
\theoremstyle{thmstyletwo}%

\theoremstyle{thmstylethree}%

\newlength\figureheight
\newlength\figurewidth 

\ifpdf
\DeclareGraphicsExtensions{.eps,.pdf,.png,.jpg}
\else
\DeclareGraphicsExtensions{.eps}
\fi

\def\norm#1{\|#1\|} 

\definecolor{darkgreen}{rgb}{0.0,0.5,0.0}
\definecolor{amber}{rgb}{1.0, 0.75, 0.0}

\def\norm#1{\left\|#1\right\|} 

\newcommand{\R}{\ensuremath{\mathbb{R}}}

\newcommand{\Ktb}{\tilde{\mathbf{K}}}
\newcommand{\Ktbt}{\tilde{\mathbf{K}}(\theta)}
\newcommand{\Ab}{\mathbf{A}}

\newcommand{\bb}{{\bf b}}

\def\matlab{{\large {\sc matlab}}}

\def\norm#1{\left\|#1\right\|} 

\DeclareMathOperator*{\argmin}{arg\,min}

\numberwithin{theorem}{section}

\begin{document}

\title[Efficient training of Gaussian processes with tensor product structure]{Efficient training of Gaussian processes with tensor product structure}


\author[1]{\fnm{Josie} \sur{Koenig}}\email{josie.koenig@uni-potsdam.de}

\author*[2]{\fnm{Max} \sur{Pfeffer}}\email{m.pfeffer@math.uni-goettingen.de}

\author[3]{\fnm{Martin} \sur{Stoll}}\email{martin.stoll@mathematik.tu-chemnitz.de}

\affil[1]{\orgdiv{Department of Mathematics}, \orgname{University of Potsdam}, \orgaddress{Karl-Liebknecht-Str. 24-25, 14476 Potsdam, \country{Germany}}}

\affil*[2]{\orgdiv{Institute of Numerical and Applied Mathematics}, \orgname{University of G\"ottingen}, \orgaddress{Lotzestr. 16-18
37083 G\"ottingen, \country{Germany}}}

\affil[3]{\orgdiv{Department of Mathematics}, \orgname{TU Chemnitz}, \orgaddress{Reichenhainer Str. 41, 09126 Chemnitz, \country{Germany}}}


\abstract{To determine the optimal set of hyperparameters of a Gaussian process based on a large number of training data, both a linear system and a trace estimation problem must be solved. In this paper, we focus on establishing numerical methods for the case where the covariance matrix is given as the sum of possibly multiple Kronecker products, i.e., can be identified as a tensor. As such, we will represent this operator and the training data in the tensor train format. Based on the \textsc{AMEn} method and Krylov subspace methods, we derive an efficient scheme for computing the matrix functions required for evaluating the gradient and the objective function in hyperparameter optimization.}

\keywords{Gaussian process, Tensor train, Trace estimation}



\maketitle

\section{Introduction}
Gaussian processes are a well-established method in machine learning and statistics to solve regression or classification problems \cite{RaWi,RaWi2}. Their performance crucially depends on the choice of the hyperparameters of the covariance kernel function and their optimization is computationally expensive and requires advanced tools from large-scale numerical linear algebra. We here focus on the case where the kernel function has additional structure and as such leads to structured covariance matrices of very large dimension. One such example given in \cite{alvarez2012kernels} are multi-output kernel functions. We discuss the efficient training if the structure of the kernel function and thus the covariance matrix or its approximation can be viewed as a tensor \cite{KoBa}. Here, we consider a tensor to be a multidimensional generalization of a matrix, i.e., an array with $D$ indices:
\begin{align}\label{Tensor}
\mathbf{K} \in \mathbb{R}^{n_1 \times \ldots \times n_D}.
\end{align}
Since the storage requirements of the tensor $\mathbf{K}$ depend on the mode $D$ in an exponential fashion, one typically observes the \textit{curse of dimensionality.} As a result, the tensor will be approximated using a low-rank tensor format \cite{Os,KoBa}. The parameter optimization requires the solution of linear systems, the approximation of a log-determinant, often reformulated as a trace estimation problem, and the computation of the gradient that then again requires linear system solves and trace estimation procedures.

Our paper starts by recalling some of the basics of Gaussian processes in Sec.~\ref{sec:GPlearning}, where also the Kronecker-sum structure of the covariance matrix will be introduced. In Sec.~\ref{sec::tensors}, we recall the basics of the tensor train (TT) format. Sec.~\ref{sec:TTKryl} describes the Krylov method in TT format. Finally, in Sec.~\ref{sec:minloglik}, we derive formulas that are necessary to compute the cost function of the hyperparameter optimization as well as its gradient, before we show numerical experiments in Sec.~\ref{sec:NumExp}.

\paragraph{Existing work}
For linear systems with the kernel matrix $\mathbf{K}$, i.e. $D=2$, the recent survey \cite{stoll2020literature} contains a list of references on this topic. One of the main ingredients is the acceleration of the matrix vector products with general kernel matrices $\mathbf{K}$. For these matrices much research is devoted to using low-rank approximations \cite{alaoui2015fast,cai2022fast,park2022randomized,rahimi2007random,wilson2015kernel}, methods from Fourier analysis \cite{nestler2021learning}, or hierarchical matrices \cite{iske2017hierarchical} in general kernel-based learning. 

In more detail, low-rank techniques have been used very successfully in Gaussian process methods often based on a set of inducing points with the \textit{subset of regressors} (SOR) \cite{silverman1985some} or its diagonal correction, the \textit{fully independent training conditional} (FITC) \cite{snelson2005sparse}. Wilson and Nickisch introduce a technique based on kernel interpolation in \cite{wilson2015kernel}. The authors in \cite{Dong} exploit this structured kernel interpolation for approximating the linear system solves as well as providing trace estimators for functions of $\mathbf{K}$ based on Krylov subspace methods. A more sophisticated implementation based on PyTorch was given in \cite{gardner2018gpytorch}. 
The challenging problem in computational Gaussian process learning and its parameter tuning is the trace estimation. This task has received much attention recently \cite{Skor,ubaru2017fast,meyer2021hutch++,cortinovis2021randomized,persson2022improved}, where the finite approximation of the expectation of $x^TAx$ is approximated. Recently, variance-reducing techniques for this approximation have been introduced to give the \textit{Hutch++} trace estimator or via the use of a preconditioner (cf. \cite{wenger2022preconditioning}). 

In the case of the underlying structure being based on low-rank tensor approximations only a few results are available, i.e., approximating the eigenfunctions of the kernel via tensor networks \cite{menzen2023projecting}, approximating the kernel matrix using a rank-one Kronecker product \cite{yu2018tensor} and for inducing point approximations of the covariance matrix \cite{izmailov2018scalable}.

\section{Gaussian process learning}
\label{sec:GPlearning}

The default setting for Gaussian processes is as follows: A set of training inputs $\mathbf X_{train} \in \mathbb R^{N \times n_{train}}$ with the corresponding outputs $\mathbf y_{train} \in \mathbb R^{n_{train}}$ is given. We assume the outputs to be polluted by centered Gaussian noise, i.e., $\mathbf y_{train}=\mathbf y_{true}+\bm{\varepsilon}$, where $\bm{\varepsilon}\sim \mathcal{N}(\mathbf{0},\sigma^2 \mathbf{I})$. The relationship between inputs and outputs must be inferred from the structure of the data. Finally, the goal is to predict the outcomes $\mathbf y_{test} \in \mathbb R^{n_{test}}$ for unseen test inputs $\mathbf X_{test} \in \mathbb R^{N \times n_{test}}$ as accurately as possible. In this paper we assume that the input-output relationship can be modeled by a Gaussian process.

Gaussian processes are a special subclass of stochastic processes. For $\mathcal{X}$ a parameter space (typically \,$\mathbb{R}^N$) and $(\Omega, \Sigma, p_{\Omega})$ a probability space, a \textit{stochastic process} is a family $\{y(x), x \in \mathcal{X}\}$ of random variables defined on $\Omega$ \cite[p.\@ 11]{LiRoSa}. This family can equivalently be described as a family of mappings $y: \mathcal{X} \times \Omega \mapsto \mathbb{R}, (x,\omega)\rightarrow y(x,\omega)$. For fixed $\omega \in \Omega$ we obtain functions on $\mathcal{X}$. For fixed $x \in \mathcal{X}$ we obtain random variables $y(x): \Omega \mapsto \mathbb{R}$. Hence, the terms \textit{random field} or \textit{random function} are synonymously used for stochastic processes.

A \textit{Gaussian process} is a stochastic process where for each finite subset $\textbf{X} \subset \mathcal{X}$, the points $(y(x))_{x \in \textbf{X}} \subset \mathbb{R}^{|\textbf{X}|}$ have a joint Gaussian distribution. Since Gaussian distributions are completely defined by their mean and variance, Gaussian processes are completely defined by their mean function $m(x)$ and their covariance or kernel function $k(x,x^\prime)$, where
\begin{align*}
m(x)&:= \mathbb{E}[(y(x))],\notag\\
k(x,x^\prime)&:= \mathbb{E}[(y(x)-m(x))(y(x^\prime)-m(x^\prime))].
\end{align*}
Since we have assumed that the input-output relationship in our data can be modeled by a Gaussian process, we now need to find the Gaussian process, i.e., the mean function and the covariance function, that best describe the training data. Of course, it is not feasible to search among all possible Gaussian processes, but we have to make some restrictions to actually get a solution. Thus, without loss of generality, it is usually assumed that $m(x) = 0 \enspace \forall x \in \mathcal{X}$ (otherwise shift by the mean function). Only the covariance function $k(x,x^\prime)$ remains as an interesting property of a Gaussian process. When it comes to kernel functions, we usually restrict ourselves to a single family of functions. We then look for the best kernel function in that family.

There are many possible families of covariance or kernel functions $k: \mathbb{R}^{\mathrm{N}} \times \mathbb{R}^{\mathrm{N}} \rightarrow \mathbb{R}$, depending on the model of choice, for example: The squared-exponential covariance functions, also called radial basis function (RBF) kernels,
\begin{equation}
k_{SE}(x,x^\prime) = \sigma_f^2\exp\Big(-\frac{\|x-x^\prime\|_2^2}{2\ell^2}\Big),
\end{equation}
with length scale $\ell >0$ and signal variance $\sigma_f>0$ as hyperparameters or the linear covariance function
\begin{equation*}
k_{lin}(x,x^\prime) = \sum_{d=1}^N \sigma^2_d x_d x_d^\prime,
\end{equation*}
with scale parameters $\{\sigma_d\}_{d=1}^N$. Here, $x_d$ denotes the $d$-th coordinate of $x$. The above examples show that the kernel functions depend on hyperparameters that can be configured to fit the data.

\paragraph{The learning procedure}
Gaussian process learning means finding the optimal hyperparameters of a fixed family of kernel functions $k_{\theta}(\cdot, \cdot)$. Here we write the parameters of the kernel to be collected in the parameter vector $\theta$. Note that we use $\tilde{\mathbf{K}}:=\mathbf{K}+\sigma^2 \mathbf{I}$ to include the noise and that the matrix $\mathbf{K}$ is obtained by evaluating $k_{\theta}(\cdot, \cdot)$ on each pair of training data. The kernel $\mathbf{K}$ depends on the parameters $\theta,$ which we write explicitly in situations where it is crucial.

The first step in the Gaussian process learning procedure is to choose a family of kernel functions to consider (e.g., squared-exponential). The goal is to find hyperparameters ${\theta^*}$ such that the Gaussian distribution ${p(\mathbf y_{train}|\mathbf X_{train})=\mathcal{N}(\bm{0},\tilde{\mathbf{K}}(\theta)) =: Z(\theta)}$ has maximum likelihood for ${\theta=\theta^*}$. The distribution $p(\mathbf y_{train}|\mathbf X_{train})$ is derived from the Gaussian process describing the input-output relationship in the training data.

Instead of maximizing the Gaussian likelihood directly, we minimize the negative logarithm of the likelihood function, called the \textit{negative log-likelihood}:
\begin{align}\label{eq::NLL}
f (\theta)= - \log Z(\theta) = \frac{1}{2} \bigg(n_{train}\log(2\pi) + \underbrace{\log\det(\tilde{\mathbf{K}})(\theta)}_{\substack{=\,\mathrm{tr}(\log\tilde{\mathbf{K}})(\theta)}} +\textbf{y}_{train}^{\mathsf{T}}\tilde{\mathbf{K}}(\theta)^{-1}\textbf{y}_{train} \bigg).
\end{align}
In order to minimize the negative log-likelihood function, we must be able to evaluate it frequently and fast. The cost of this evaluation is dominated by solving a linear system $\tilde{\mathbf{K}}^{-1}(\theta)\textbf{y}_{train}$ and computing its log-determinant $\log\det(\tilde{\mathbf{K}}(\theta)) = \mathrm{tr}(\log\tilde{\mathbf{K}}(\theta))$. To minimize \eqref{eq::NLL} the gradients of the function $f(\theta)$ are also of interest. It is well known that the analytical gradients of the log-determinant and the kernel inverse are given by 
\begin{equation}
\label{eq::derK}
\frac{\partial}{\partial \theta_j}\log\det(\Ktbt)=\mathrm{tr}\left(\Ktbt^{-1}\frac{\partial \Ktbt}{\partial \theta_j}\right),\
\frac{\partial}{\partial \theta_j}\Ktbt^{-1}=-\Ktbt^{-1}\frac{\partial \Ktbt}{\partial \theta_j}\Ktbt^{-1}.
\end{equation}
These computations necessary for the minimization are typically difficult to handle, especially when the amount of training data is large. Efficient optimization of the parameters with respect to the training data will be the focus of our investigation.

\paragraph{Predictions from unseen data}
After optimizing the hyperparameters, we obtain the kernel function $k_{\theta^*}(\cdot, \cdot)$ that best describes our training data.
We can then obtain predictions $\mathbf y_{test}$ from new unseen input data $\mathbf X_{test}$ using the learned kernel in the following way:


The prior joint distribution of the noise-polluted outputs $\mathbf y_{train}$ and $\mathbf y_{test}$ using the optimized kernel is again Gaussian:
\begin{align*}
\begin{bmatrix}
\mathbf y_{train}\\
\mathbf y_{test}
\end{bmatrix}
\sim
\mathcal{N}
\left(
\mathbf{0},
\begin{bmatrix}
\tilde{\mathbf{K}}&\mathbf{K}_*\\
\mathbf{K}_*^\top&\mathbf{K}_{**}+\sigma^2 \mathbf{I}
\end{bmatrix}
\right),
\end{align*}
where $\tilde{\mathbf{K}}$ is as above and $\mathbf{K}_*=k_{\theta^*}(\mathbf X_{test},\mathbf X_{train})$ and $\mathbf{K}_{**}=k_{\theta^*}(\mathbf X_{test},\mathbf X_{test})$ are the matrices obtained by evaluating the learned kernel function on the data.
The predictive posterior distribution for the noisy outputs $\mathbf y_{test}$ from the test points $\mathbf X_{test}$ can then be computed by
\begin{align*}
p(\mathbf y_{test}|\mathbf X_{test}, \mathbf X_{train},\mathbf y_{train})&=\mathcal{N}(\mathbf y_{test}|\bm{\mu}_*,\mathbf{\Sigma}_*),\\
\bm{\mu}_*&=\mathbf{K}_*^\top (\mathbf{K}+\sigma^2 \bm{I})^{-1}\mathbf y_{train}\,\\
\bm{\Sigma}_*&=\mathbf{K}_{**}+\sigma^2 \mathbf{I}-\mathbf{K}_*^\top (\mathbf{K}+\sigma^2 \bm{I})^{-1}\mathbf{K}_*.
\end{align*}
Thus, one needs to compute the values of the learned kernel function for each pair of inputs to obtain the mean $\bm{\mu}_*$ and the covariance $\bm{\Sigma}_*$ of the test outputs. The desired prediction $\mathbf y_{test}$ is then obtained by sampling from the posterior Gaussian distribution $\mathcal{N}(\bm{\mu}_*,\mathbf{\Sigma}_*)$. Alternatively, the mean $\bm{\mu}_*$ can be used as a point estimate for $\mathbf y_{test}$ and the covariance $\mathbf{\Sigma}_*$ describes its uncertainty. {Fig.~\ref{fig:gp} illustrates this prediction procedure for a one-dimensional example. Fig.~\ref{subfig-1:GP} shows the prior distribution, Fig.~\ref{subfig-2:GP} shows the posterior distribution after the training data has been incorporated.

\begin{figure}[t]
\subfloat[Gaussian process prior\label{subfig-1:GP}]{%
\includegraphics[width=0.5\textwidth]{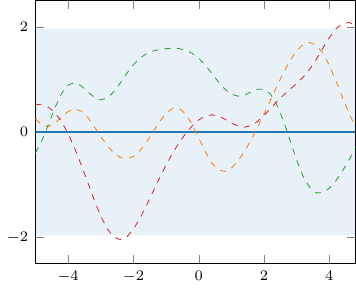}
}
\hfill
\subfloat[Gaussian process posterior with small noise\label{subfig-2:GP}]{%
\includegraphics[width=0.5\textwidth]{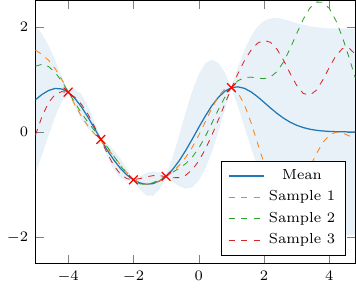}
}
\caption{(a) Gaussian process prior with no training data. (b) The posterior distribution after having incorporated $5$ training points. Means are given in blue, the shaded area indicates the $95\%$ confidence interval and the dashed lines sample from the current distribution.}
\label{fig:gp}
\end{figure}

The above predictive procedure requires knowledge of the optimal kernel function $k_{\theta^*}(\cdot, \cdot)$, and the quality of the prediction depends crucially on its hyperparameters. From Fig.~\ref{fig:gppara} it can be seen how the choice of parameters affects the posterior distribution and the ability to make meaningful predictions for unseen data. This emphasizes why the hyperparameter optimization described above is essential.

\begin{figure}[t]
\subfloat[Small length scale\label{subfig-1:GPpara}]{%
\includegraphics[width=0.5\textwidth]{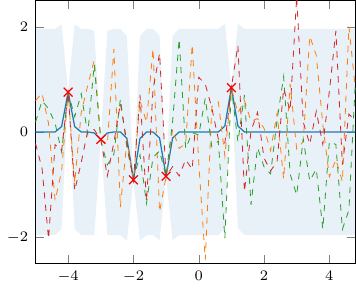}
}
\hfill
\subfloat[Large length scale\label{subfig-2:GPpara}]{%
\includegraphics[width=0.5\textwidth]{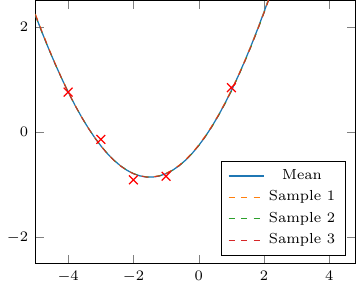}
}

\subfloat[Small signal variance\label{subfig-3:GPpara}]{%
\includegraphics[width=0.5\textwidth]{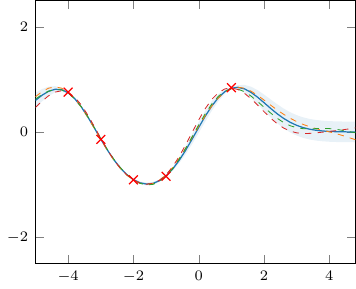}
}
\hspace*{-.2cm}
\subfloat[Large signal variance\label{subfig-4:GPpara}]{%
\includegraphics[width=0.51\textwidth]{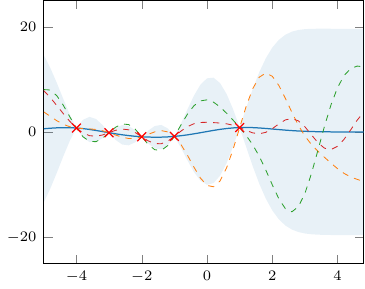}
}
\caption{Gaussian posterior from a squared-exponential kernel function with different length scale and signal variance. Training and test data are the same, given by red crosses. Means are given in blue, the shaded areas indicates the $95\%$ confidence interval and the dashed lines sample from the current distribution. Note the different y-scale in (d).}
\label{fig:gppara}
\end{figure}

\paragraph{The Kronecker-sum kernel}
The numerical methods to be used in hyperparameter optimization depend on the chosen kernel. In this paper, we consider a special setting for input locations on a Cartesian grid, cf.\@ \cite[p.\@ 126 ff.]{Saa}. The test and training data must be of the form
\begin{align*}
\textbf{X} = \textbf{X}^{(1)} \times \ldots \times \textbf{X}^{(D)}, \enspace \textbf{X}^{(d)}\in \mathbb{R}^{n_d},\enspace d = 1,\ldots,D,
\end{align*}
where $\textbf{X}^{(d)}$ contains the $n_d$ input locations along the dimension $d$ which may vary for each $d = 1,\ldots,D$. In total, there are $n = \prod_{d=1}^D n_d$ inputs for $\textbf{X} \subset \mathbb{R}^D$ leading to outputs $\mathbf{y}\in \mathbb{R}^n = \mathbb{R}^{n_1 \times \cdots \times n_D}$. Outputs of this form can be treated as $D$-mode tensors and will be used as such in our computations.  

With this in mind, we can now define the covariance kernel with the following Kronecker-sum structure: 
\begin{align}\label{eq::kronCov}
\mathbf{K}(\mathbf{X},\mathbf{X})&=\sum_{r=1}^{R}\mathbf{K}_r^{(1)}\otimes\mathbf{K}_r^{(2)}\otimes\ldots \otimes\mathbf{K}_r^{(D)}, \text{ where }  \mathbf{K}^{(d)}_r:= k^{(d)}_r(\textbf{X}^{(d)},\textbf{X}^{(d)}).
\end{align}
Such a kernel corresponds to a $2D$-mode tensor with each $\mathbf{K}_r^{(d)}\in\R^{n_d\times n_d}.$ In particular, we will focus on the case when $R>1$, since the case $R=1$ with moderate individual dimensions of the kernel matrices $\mathbf{K}_1^{(d)}$ can be treated with Cholesky decompositions of these matrices \cite{saad2003iterative}. In that case, efficient training of the Gaussian process could be performed using standard techniques.

The ideas presented in this paper apply in principle to all kinds of kernel functions within the Kronecker-sum kernel given in \eqref{eq::kronCov}. However, from now on we will assume that the one-dimensional covariance matrices $\mathbf{K}^{(d)}_r$ result from applying a squared-exponential covariance function ${k^{(d)}_r: \mathbb{R}\times\mathbb{R}\rightarrow\mathbb{R}}$ with length scale parameter $\ell^{(d)}_r$ to each pair of inputs $x,x^\prime \in \textbf{X}^{(d)}$. The signal variance $\sigma_{f(r)}$ can be shifted into the first kernel function $k^{(1)}_r$, so we use:
\begin{align}\label{eq::single_covs}
\begin{split}
k_r^{(1)}(x,x^\prime) &= \sigma_{f(r)}^2\exp\Big(-\frac{(x-x^\prime)^2}{2(\ell^{(1)}_r)^2}\Big),\quad r = 1,\ldots,R, \enspace \text{ and }\\
k_r^{(d)}(x,x^\prime) &= \exp\Big(-\frac{(x-x^\prime)^2}{2(\ell^{(d)}_r)^2}\Big), \quad r = 1,\ldots,R,\enspace d = 2,\ldots,D.
\end{split}
\end{align}
For the setting described above we need to learn the hyperparameters collected in ${\theta = \{\sigma_{f(r)},\ell^{(d)}_r\}_{r = 1,\ldots,R,\, d = 1,\ldots,D}}$ and each hyperparameter appears only in one dimension and one summand. This will be especially important when computing the kernel derivative $\frac{\partial \Ktbt}{\partial \theta_j}$ with respect to the different hyperparameters. For the squared exponential kernel ${k_{SE}(x,x^\prime) = \sigma_f^2\exp\Big(-\frac{\|x-x^\prime\|_2^2}{2\ell^2}\Big)}$, it is advantageous to use the logarithmic transformation (\textit{log-transform}) of the hyperparameters \cite[p.\@ 138]{Saa}. It holds
	\begin{align*}
	\frac{\partial k_{SE}(x,x^\prime)}{\partial\sigma_f} = \frac{2 k_{SE}(x,x^\prime)}{\sigma_f} \enspace \text{  and  } \enspace \frac{\partial k_{SE}(x,x^\prime)}{\partial\ell} = \frac{\|x-x^\prime\|_2^2}{\ell^3}k_{SE}(x,x^\prime),
	\end{align*}
	but due to the chain rule with the derivative of the natural logarithm we get
	\begin{align}\label{eq::RBF_derivs}
	\frac{\partial k_{SE}(x,x^\prime)}{\partial\log(\sigma_f)} = 2 k_{SE}(x,x^\prime) \enspace \text{  and  } \enspace \frac{\partial k_{SE}(x,x^\prime)}{\partial\log(\ell)} = \frac{\|x-x^\prime\|_2^2}{\ell^2}k_{SE}(x,x^\prime).
	\end{align}
	The main advantages of the log-transform are that the positivity requirements for $\sigma_f$ and $\ell$ are automatically satisfied, and that dividing by $\ell^2$ instead of $\ell^3$ and not dividing by $\sigma_f$ makes the derivative computation numerically more stable, especially for hyperparameters near zero. The values of the derivatives using the log-transform will be scaled relative to the true values but the roots will remain the same, i.e., $\frac{\partial k_{SE}(x,x^\prime)}{\partial\sigma_f} = 0$ iff $\frac{\partial k_{SE}(x,x^\prime)}{\partial\log(\sigma_f)} = 0$ and the same for $\ell$. This justifies replacing the hyperparameters with their log-transform in the derivative computation.
	
Kernels with Kronecker-sum structure allow for the modeling of overlapping phenomena, as each summand models an independent stochastic process. Such settings occur, for example, in spatio-temporal magnetoencephalography (MEG) \cite{bijma2005MEG} or climate data sets \cite{Saa}. Computations involving high-dimensional Kronecker-sum kernels require advanced numerical algebra techniques, which will be discussed in the remainder of this paper.

\section{A low-rank tensor train framework}
\label{sec::tensors}
Since our covariance operator is given in tensor product form we can use low-rank tensor formats to perform the necessary operations. Computation and storage of higher order tensors suffer from the {\em curse of dimensionality} since their complexity grows exponentially with the order of the tensor. Similar to~\eqref{eq::kronCov}, a tensor $\mathbf T \in \mathbb R^{n_1 \times \cdots \times n_D}$ can be expressed as a sum of elementary tensors
\begin{equation}
\mathbf T = \sum_{r=1}^R \mathbf t_1^{(r)} \otimes \cdots \otimes \mathbf t_D^{(r)},
\end{equation}
where $\mathbf t_d^{(i)} \in \mathbb R^{n_d}$ for $d = 1,\ldots,D$ and $r = 1,\ldots,R$. If $R$ is the smallest integer for which such a representation exists, it is called the {\em tensor rank} of the tensor and the representation is called a {\em canonical decomposition} of $\mathbf T$. If this rank is small, the complexity of storage is now linear in $n = \max\{n_1,\ldots,n_D\}$.

While the canonical decomposition format is concise and efficient, it possesses a number of well-documented drawbacks that make its use in optimization methods tedious: The computation of the tensor rank (and therefore of the decomposition itself) is NP-hard and, furthermore, the set of tensors with tensor rank at most $R$ is not closed, meaning that we can find tensors of a certain tensor rank that can be approximated arbitrarily well by tensors of lower tensor rank (cf.~border rank problem \cite{landsbergtensors}). In short, if the canonical decomposition of a tensor is unknown, the attempt to compute it will often result in failure.

For this reason, other tensor decomposition formats have been proposed in the literature. The most well-suited class of decomposition formats for the purpose of optimization are so called {\em Tree Tensor Networks} \cite{Hackbusch2009, Os}. These can be computed efficiently using higher order generalizations of the SVD, quasi-best low-rank approximations are readily available \cite{Lathauwer2000}, and tensors with given multilinear rank form smooth manifolds that are embedded in the tensor ambient space \cite{Holtz2012b, USCHMAJEW2013133}. In this article, we will focus on the tensor train (TT) decomposition, as it is a good compromise between simplicity (of the format) and complexity (of computation and storage).

The TT decomposition of an order-$D$ tensor $\mathbf T \in \mathbb R^{n_1 \times \cdots \times n_D}$ can be given elementwise as 
\begin{align}\label{TTStand}
\textbf{T}(i_1,\ldots,i_D) = \sum_{k_1=1}^{r_1}\cdots\sum_{k_{D-1}=1}^{r_{D-1}}\textbf{T}_1(i_1,k_1)\textbf{T}_2(k_1,i_2,k_2)\cdots \textbf{T}_D(k_{D-1},i_{D-1}),
\end{align}
where $\textbf{T}_{d} \in \mathbb{R}^{r_{d-1}\times n_d \times r_d}$ for $d = 1,\ldots,D$ and $r_0 = r_D = 1$. The tuple ${\mathbf r = (1,r_1,\ldots,r_{D-1},1)}$ is called the \textit{TT rank} of the tensor $\textbf{T}$. The tensors $\textbf{T}_d$ of order (up to) three are called the \textit{TT cores} of the TT decomposition. The complexity of storage of a TT tensor is therefore no longer exponential in the order $D$ but linear in $D$ and $n$ and quadratic in the ranks.


Following the same idea, we decompose a linear operator $\textbf{A} : \mathbb R^{n_1 \times \cdots \times n_D} \rightarrow \mathbb R^{m_1 \times \cdots \times m_D}$ in a TT-like format. For a fixed basis and after some reordering of indices, this operator can be given as a tensor of order $2D$ with mode sizes $m_1,n_1,\ldots,m_D,n_D$. Analogously to the TT format, we can therefore define a TT matrix format by
\begin{align}\label{TTOp}
\textbf{A}(i_1,j_1,\ldots,i_D,j_D) = \sum_{k_1=1}^{r_1}\cdots\sum_{k_{D}=1}^{r_{D}}\bigotimes_{d=1}^D \textbf{A}_d(k_{d-1},i_d,j_d,k_d),
\end{align}
where $\textbf{A}_d \in \mathbb{R}^{r_{d-1}\times m_d \times n_d \times r_d}$ for $d = 1,\ldots,D$ and again $r_0 = r_D = 1$.

For two TT tensors of same order and dimension, but with different TT ranks, most standard operations can be computed efficiently. This includes addition of two tensors, the Hadamard product, and the Frobenius inner product between them (denoted by $\mathtt{dot}$). Furthermore, given a TT matrix of suitable column dimension, we can efficiently compute the matrix-vector product with a TT tensor by computing matrix-vector products of the cores. These operations are implemented in {\matlab}, e.g., in the $\mathtt{tt-toolbox}$~\cite{TTTB}. For the exact definitions of the above operations, we refer the reader to~\cite{Os}.

The complexity of computations and storage is only efficiently reduced if the TT ranks are low. Most of the above operations increase the ranks, sometimes drastically. Therefore, we employ a rounding procedure, known as {\em TT-SVD}, that truncates the ranks when they become exceedingly large. The TT-SVD effectively relies on successive truncated SVDs of the reshaped TT cores. One sweep through all cores reduces the ranks to a desired magnitude or until a designated error threshold has been reached. We denote this rounding procedure by $\mathtt{round}$ and we predefine the desired error tolerance.

Many algorithms prove to be stable with respect to this rank rounding procedure and the errors remain moderate. Furthermore, given a fixed TT rank ${\mathbf r = (1,r_1,\ldots,r_{D-1},1)}$, the TT-SVD produces a quasi-best rank-$\mathbf r$ approximation of $\mathbf T$, meaning that the error differs only by a factor $\sqrt D$ from the error of the best rank-$\mathbf r$ approximation of $\mathbf T$. 

For the implementation, we use the $\mathtt{round}$ procedure from the $\mathtt{tt-toolbox}$ in {\matlab} \cite{TTTB}. Furthermore, from the same toolbox, we use the \textsc{AMEn}-method for the solution of a linear system in TT format with rank control (denoted by $\mathtt{amen}$) and also for the summation of a collection of TT tensors (denoted by $\mathtt{amen\_sum}$).

Since the training data $\mathbf{y}_{train} \in \mathbb R^{n_1 \times \cdots \times n_D}$ is a tensor, we can bring it into TT format using the (possibly truncated) SVD. Each summand $\mathbf{K}_r^{(1)}\otimes\mathbf{K}_r^{(2)}\otimes\ldots \otimes\mathbf{K}_r^{(D)}$ of the kernel matrix $\mathbf{K}$ is a rank-$1$ TT matrix and we can use summation of TT matrices to obtain the kernel in TT format (usually without rounding).



\section{The TT-Krylov method}\label{sec:TTKryl}


In this section, we introduce the {\em TT-Krylov method} for the efficient computation of $g(\Ab)\mathbf{b}$ and $\mathbf{v}^Tg(\Ab)\mathbf{v}$ for (possibly nonsymmetric) TT matrices $\Ab$, TT tensors $\mathbf b, \mathbf v$, and $g$ a matrix function. 

We start with a general introduction for the matrix case. Even then,
this problem poses a significant challenge as the evaluation of matrix functions is expensive, especially once the matrix dimensions are large \cite{higham2008functions,guttel2013rational}. Many efficient techniques exist for evaluating a matrix function times a vector, i.e., $g(\Ab)\mathbf{b}$. Problems of the form $\textbf{v}^Tg(\Ab)\textbf{u}$ (often with $\mathbf{u}=\mathbf{v}$) have been studied in the seminal works of Golub and Meurant \cite{golub2009matrices,golub2020matrices,golub1997matrices} that are rooted in the close connection of this expression to Gaussian quadrature. 

In this paper we focus on techniques based on Krylov subspaces and briefly illustrate our approach here. The Arnoldi method is used to generate an orthonormal basis for the Krylov subspace
$$
\mathcal{K}_\ell(\Ab,\mathbf{b})=\mathrm{span}\left\lbrace \mathbf{b},\Ab\mathbf{b},\ldots, \Ab^{\ell-1}\mathbf{b}\right\rbrace
$$
via a Gram-Schmidt procedure. One obtains the following decomposition
$$
\Ab \mathbf{V}_\ell=\mathbf{V}_\ell\mathbf{H}_{\ell}+h_{\ell+1,\ell} \mathbf{v}_{\ell+1}\mathbf{e}_\ell^{\mathsf T},
$$
where the vectors $\mathbf{v}_k$ are the basis vectors of the Krylov subspace, $\mathbf{V}_\ell=[\mathbf v_{1} \cdots \mathbf v_{\ell}],$ and $\mathbf{H}_{\ell}\in\R^{\ell \times \ell}$ a Hessenberg matrix. With the choice of selecting the seed vector for the Krylov subspace to be 
$$
\mathbf v_1=\bb/\norm{\bb}
$$
we obtain the following approximation
$$
g(\Ab)\mathbf{b}\approx \norm{\bb}\mathbf{V}_\ell g(\mathbf{H}_{\ell})\mathbf{e}_1,
$$
which is based on the approximation resulting from the Arnoldi procedure.
Since the dimensionality of $\mathbf{H}_{\ell}$ is much smaller than the dimensionality of $\mathbf{A}$, $g(\mathbf{H}_{\ell})$ can now be evaluated efficiently \cite{higham2008functions}. 

An alternative to standard Krylov methods are rational Krylov methods \cite{guttel2013rational,guttel2010rational,ruhe1984rational}. These methods approximate the computation of $g(\mathbf{A})\mathbf{b}$ via a rational polynomial $r_\ell(\mathbf{A})$ via $r_\ell(\mathbf{A})\mathbf{b},$ where $r_\ell=\frac{p_{\ell-1}}{q_{\ell-1}}$ and $q_{\ell-1}(z)=\prod_{j=1}^{\ell-1}(1-\frac{z}{\xi_j})$ with poles $\xi_j$. The rational Krylov space is then given via 
$
\mathcal{Q}_\ell(\Ab,\mathbf{b})=q_{\ell-1}(\Ab)^{-1}\mathrm{span}\left\lbrace \mathbf{b},\Ab\mathbf{b},\ldots, \Ab^{\ell-1}\mathbf{b}\right\rbrace.
$
The corresponding Arnoldi method then reads in matrix form as
$
\mathbf{A}\mathbf{V}_{\ell+1}\mathbf{K}_\ell=\mathbf{V}_{\ell+1}\mathbf{H}_\ell
$
where $\mathbf{K}_\ell,\mathbf{H}_\ell$ are $\ell+1$ by $\ell$ unreduced Hessenberg matrices. The approximation is then obtained via
$
g(\mathbf{A})\mathbf{b}\approx
\mathbf{V}_\ell g(\mathbf{A}_\ell)\mathbf{V}_\ell^{T}\mathbf{b},
$
where the Rayleigh quotient $\mathbf{A}_\ell=\mathbf{V}_\ell^{T}\mathbf{AV}_\ell$ can also be written as
$\mathbf{A}_\ell=\mathbf{H}_\ell\mathbf{K}_\ell^{-1}.$ As a result we can get the approximation from the fact that 
$
g(\mathbf{A})\mathbf{b}\approx \norm{\mathbf{b}}
\mathbf{V}_\ell g(\mathbf{A}_\ell)\mathbf{e}_1.
$

In particular, if we are interested in an computing the value of the term $\mathbf{v}^{\mathsf T}g(\Ab)\mathbf{v}$  we get as an approximation from the Arnoldi procedure  
$$
\mathbf{v}^{\mathsf T}g(\Ab)\mathbf{v}\approx \norm{\mathbf v}^2\mathbf{e}_1^{\mathsf T} g(\mathbf{H}_{\ell})\mathbf{e}_1,
$$
and from the rational Arnoldi we see
$$
\mathbf{v}^{\mathsf T}g(\Ab)\mathbf{v}\approx \norm{\mathbf v}^2\mathbf{e}_1^{\mathsf T} g(\mathbf{A}_{\ell})\mathbf{e}_1.
$$

One can also rely on approximations based on the nonsymmetric Lanczos process \cite{strakovs2011efficient,schweitzer2017two} but these are more delicate for the evaluation of expressions of the form $\mathbf{u}^{\mathsf T}g(\Ab)\mathbf{v}$ as one typically requires $\mathbf{u}^{\mathsf T}\mathbf{v}\neq 0,$ which is not satisfied in our case (see \eqref{approxFAB}).

For the approximation of the gradient and trace estimation we employ the Krylov and rational Krylov method implemented in the tensor train format. For this, the steps of any Krylov method are performed using the tensor train arithmetic
$$
\mathbf{K} \mathbf{V}_\ell=\mathbf{V}_\ell\mathbf{H}_{\ell}+h_{\ell+1,\ell}v_{\ell+1}\mathbf{e}_\ell^{\mathsf T},
$$
where now $\mathbf{V}_\ell$ represents a collection of $\ell$ TT tensors. The matrix $\mathbf{H}_{\ell}$ here contains the coefficients resulting from the orthogonalization process within the Arnoldi method. Even though the elements in $\mathbf{V}_\ell$ are TT tensors, $\mathbf{H}_{\ell}$ remains an $\ell \times \ell$ matrix to which we apply the matrix function. This can be done in the same way for the rational Krylov methods.

However, we need to make sure that the TT ranks do not grow too large. Therefore, we round off the TT ranks after each orthogonalization step. This is done using the $\mathtt{round}$ procedure with a predefined truncation tolerance. {This affects the orthogonality of the basis vectors $\mathbf v_k$ and we therefore repeat the orthogonalization (including the rounding) once more. This is akin to a reorthogonalization step in the matrix case.} Our experiments show that then, the ranks remain managable and the error is small. Ultimately, this rounding step is the only alteration of the Krylov method in the matrix case. As a stopping criterion we use the difference between two consecutive approximations to the desired quantities. 

In the case of approximating 
$$
\mathbf{v}^{\mathsf T}g(\Ab)\mathbf{u}
$$
techniques based on the nonsymmetric Lanczos process can be applied \cite{golub2009matrices} but since these might suffer from serious breakdowns we will use the approximation $g(\Ab)\mathbf{u}$ and then compute the inner product with $\mathbf{v}$ afterwards. 

The procedure is summarized in Alg.~\ref{alg:fAb}, where we show the case for a nonsymmetric TT-matrix $\Ab$. In the symmetric case, lines $3$ to $6$ of the Algorithm would only orthogonalize for $j=k-1,k$. If we are interested in the approximation of $\mathbf{v}^{\mathsf T}g(\Ab)\mathbf{v}$ we would replace line $12$ by
$$
\mathbf{v}^{\mathsf T}g(\Ab)\mathbf{v}\approx \norm{\mathbf v}^2\mathbf{e}_1^{\mathsf T} g(\mathbf{H}_{\ell})\mathbf{e}_1,
$$
and denote the call by $\mathtt{tt\_krylov(g,\Ab,\mathbf v,\mathbf v,\mathrm{trunctol})}$.

\begin{algorithm}[H]
	\caption{$\mathtt{tt\_krylov(g,\Ab,\bb,\mathrm{trunctol})}$: Algorithm for approximating $g(\Ab)\mathbf{b}$ in TT format. The algorithm shows Arnoldi iteration for the basis creation.}\label{alg:fAb}
	\begin{algorithmic}[1]
		\Require $g,\Ab,\bb,\mathrm{trunctol}$
		\Ensure $g(\Ab)\bb$
		\For {$k=1:\mathrm{maxit}$}
		\State $ \mathbf{w}_{k}=\Ab \mathbf{v}_{k}$
		\For {two loops}
		\For {$j=1:k$}
		\State $h_{j,k}= h_{j,k}+ \mathtt{dot}(\mathbf{v}_{j},\mathbf{w}_{k})$
		\State $\mathbf{w}_{k} = \mathtt{round}(\mathbf{w}_{k}-h_{k,j}\mathbf{v}_{j},\mathrm{trunctol})$
		\EndFor 
		\EndFor
		\State $h_{k+1,k}= \mathtt{dot}(\mathbf{w}_{k},\mathbf{w}_{k})^{1/2}$
		\State $\mathbf{v}_{k+1}=\frac{\mathbf{w}_{k}}{h_{k+1,k}}$  
		\State $\mathbf{gH}= g(\mathbf{H}(1:k, 1:k))$
		\State $
		g(\Ab)\mathbf{b}\approx\mathtt{amen\_sum}(\mathbf{v}_{1:k}, \mathbf{gH}_{:,1}\norm{\bb}, \mathrm{trunctol})$
		\State Check difference between to consecutive approximations and if small \textbf{stop.}
		\EndFor
	\end{algorithmic}
\end{algorithm}

\section{Minimizing the negative log-likelihood in TT format}
\label{sec:minloglik}

We recall that the task in Gaussian process learning is to minimize the negative log-likelihood~\eqref{eq::NLL} for the parameters $\theta$:
\begin{align}
\theta^* = \argmin_{\theta} \frac{1}{2} \bigg(n_{train}\log(2\pi) + \log\det(\tilde{\mathbf{K}})(\theta) + \textbf{y}_{train}^{\mathsf{T}}\tilde{\mathbf{K}}(\theta)^{-1}\textbf{y}_{train} \bigg).
\end{align}
Since the covariance matrix $\mathbf K(\theta)$ and the training data $\mathbf{y}_{train}$ can be represented in TT format, all operations in the optimization problem are performed in TT format also. As discussed, this results in a significant reduction of complexity and it allows for the solution of very large problems. For the minimization of the negative log-likelihood, we employ a basic LBFGS-solver provided by {\matlab}. For this, all we need is the efficient computation of the cost function and the gradient. This requires the evaluation of the log-determinant and the linear system (in the cost function) and the computation of the trace and another linear system (in the gradient). We describe our procedure in the following subsections.

\subsection{Computing the cost function}

\paragraph{Solving a linear system in TT format}

The first part of the objective function we discuss is the solution of the linear system that comes from the evaluation of the term $\textbf{y}^{\mathsf{T}}(\mathbf{K} + \sigma^2\textbf{I})^{-1}\textbf{y}.$ Computing the Cholesky decomposition of the matrix $\tilde{\mathbf{K}}=\mathbf{K} + \sigma^2\textbf{I}$ is too costly for large $D$. 
For this we will employ the \textsc{AMEn} method introduced in \cite{dolgov2014alternating} without any further modifications. This method efficiently approximates the solution of a linear system given in TT format, while maintaining a low TT rank. It consists of an alternating procedure that cycles through the TT components and optimizes them separately. In each step, a small fraction of the residual is used in order to adapt the TT ranks. We refer the reader to the original article for more details.

\paragraph{Trace estimation and matrix functions}

The required efficient evaluation of the log-determinant becomes intractable for large matrix dimensions and we rely on the equivalent relation
$$
\log\det(\Ktb)=\mathrm{tr}(\log(\Ktb)).
$$
Using this will allow us to avoid the computation of the determinant of $\Ktb$ altogether. Note that again, the problem would be rather easy if we could afford a Cholesky decomposition of the matrix $\Ktb$. Our goal is to estimate the trace of the matrix logarithm. For this, we use the Hutchinson trace estimator
$$
\mathrm{tr}(\Ab)=\mathbb{E}(\textbf{z}^{\mathsf T}\Ab \textbf{z})
$$
for a Rademacher or Gaussian random vector $\textbf{z}$. Naturally, we approximate this using a Monte-Carlo approach, i.e., with a finite sum
$$
\mathrm{tr}(\Ab)\approx \frac{1}{p}\sum_{i=1}^p(\textbf{z}_i^{\mathsf T}\Ab \textbf{z}_i)
$$
for $p$ i.i.d.~ Rademacher vectors $\textbf{z}_i$. For $\Ab=\log(\Ktbt)$ we can now use the trace estimator for the evaluation of the objective function. However, as is usually the case for naive Monte-Carlo style approaches, the relative error of the trace will be proportional to $\frac{1}{\sqrt{p}}$. In the literature, it is common to choose relatively low numbers of probe vectors ($p < 100$) \cite{Dong}. The error is then acceptable for the evaluation of the cost function. However, when computing the gradient (where another trace estimation will become necessary), we would need higher accuracy to ensure that the search direction is actually a descent direction of the cost function. For this reason, it is common to replace the cost function by a numerically efficient version
\begin{align}\label{eq::NLL2}
f (\theta)\approx \hat{f}(\theta):=\frac{1}{2} \left(N\log(2\pi) + \frac{1}{p}\sum_{i}(\textbf{z}_i^{\mathsf T}\log(\Ktbt) \textbf{z}_i) +\textbf{y}^{\mathsf{T}}(\Ktbt)^{-1}\textbf{y} \right).
\end{align}
In the following, we therefore minimize this modified cost function instead.

The computational challenge now lies in efficiently evaluating the quantity
$\sum_{i}(\textbf{z}_i^{\mathsf T} \log({\Ktb}) \textbf{z}_i)$. This is done using the symmetric TT-Krylov method introduced in Sec.~\ref{sec:TTKryl}. We also need to convert the probe vectors $\textbf{z}_i$ to TT format. However, these vectors will usually have full TT rank, and we therefore use TT tensors of rank 1, where each component is an i.i.d.~Rademacher vector:
$$
\textbf{z}_i = \textbf{z}_i^{(1)} \otimes \cdots \otimes \textbf{z}_i^{(D)}, \qquad \textbf{z}_i^{(d)} \sim \mathcal U([-1,1]^{n_d}), \; d = 1,\ldots,D.
$$
This keeps the TT ranks managable during the TT-Krylov procedure, and we obtain a framework for the efficient computation of the cost function. 

We summarize the computation of the cost function in Alg.~\ref{alg:cost}

\begin{algorithm}[H]
	\caption{Algorithm for the computation of the cost function.}\label{alg:cost}
	\begin{algorithmic}[1]
		\Require $\theta = \{\sigma_{f(r)},\ell^{(d)}_r\}_{r = 1,\ldots,R,\, d = 1,\ldots,D}, p, \mathrm{kryltol}, \mathrm{amentol}$
		\Ensure $\hat f(\theta)$
		\State {$t = 0$}
		\For {$i=1:p$}
		\State {$\textbf{z}_i = \textbf{z}_i^{(1)} \otimes \cdots \otimes \textbf{z}_i^{(D)}, \; \textbf{z}_i^{(d)} \sim \mathcal U([-1,1]^{n_d}), \; d = 1,\ldots,D$}
		\State {$t = t + \mathtt{tt\_krylov}(\mathtt{logm},\Ktbt,\mathbf z_i,\mathbf z_i,\mathrm{kryltol})$}
		\EndFor	
		\State {$\mathbf \alpha = \mathtt{amen(\Ktbt,\mathbf y_{train},\mathrm{amentol})}$}
		\State {$\hat f(\theta) = (1/2)(N \log(2\pi) + t/p + \mathbf y_{train}^{\mathrm T} \mathbf \alpha)$}
	\end{algorithmic}
\end{algorithm}

\subsection{Computing the gradient}

For the training of the hyperparameters we require the derivative of the approximate negative log-likelihood~\eqref{eq::NLL2} with respect to the parameters $\theta$. Therefore, we need the derivative of the inverse kernel as well as that of the log-determinant.

\paragraph{Derivative of the inverse kernel}

Using the inverse function rule, it holds that
$$
\frac{\partial \tilde{\mathbf{K}}(\theta)^{-1}}{\partial \theta_j} = -\tilde{\mathbf{K}}(\theta)^{-1}\frac{\partial \tilde{\mathbf{K}}(\theta)}{\partial \theta_j}\tilde{\mathbf{K}}(\theta)^{-1}.
$$
Since $\tilde{\mathbf{K}}(\theta) = \mathbf{K}(\theta) + \sigma^2\textbf{I}$, we therefore require the efficient evaluation of the quadratic expression 
$$
\textbf{y}_{train}^{\mathsf{T}} \frac{\partial \tilde{\mathbf{K}}^{-1}(\theta)}{\partial \theta_j}\textbf{y}_{train}=
-\textbf{y}_{train}^{\mathsf{T}}\tilde{\mathbf{K}}^{-1}(\theta)\frac{\partial \mathbf{K}(\theta)}{\partial \theta_j}
\tilde{\mathbf{K}}^{-1}(\theta)\textbf{y}_{train} =
-\alpha^{\mathsf{T}}\frac{\partial \mathbf{K}(\theta)}{\partial \theta_j}\alpha,
$$
where $\tilde{\mathbf{K}}(\theta)\alpha=\textbf{y}_{train}$ can be precomputed using again the \textsc{AMEn} method. The computation of the derivative $\frac{\partial \mathbf{K}(\theta)}{\partial \theta_j}$ in TT matrix format will be addressed below.

\paragraph{Derivative of the log-determinant}


If we now study \eqref{eq::derK} it might seem natural to assume that we approximate the derivative of the log-determinant using 
$$
\frac{\partial \log\det(\Ktbt)}{\partial \theta_j} = \mathrm{tr}\left(\Ktbt^{-1}\frac{\partial \mathbf{K}(\theta)}{\partial \theta_j}\right)
\approx \frac{1}{p}\sum_{i}\textbf{z}_i^{\mathsf T}\Ktbt^{-1}\frac{\partial \mathbf{K}(\theta)}{\partial \theta_j}\textbf{z}_i.
$$
But in our scheme we are interested in the gradient of $\hat{f}(\theta)$, which requires the derivative of the expression
$$
\textbf{z}_i^{\mathsf T} \frac{\partial \log(\Ktbt) }{\partial \theta_j}\textbf{z}_i.
$$
If $\Ktbt^{-1}$ and $\frac{\partial \mathbf{K}(\theta)}{\partial \theta_j}$ commute, the derivative of the matrix logarithm with respect to a parameter is (see \cite{haber2018notes} for a detailed derivation)
$$
\frac{\partial \log(\Ktbt)}{\partial \theta_j}=\Ktbt^{-1}\frac{\partial \mathbf{K}(\theta)}{\partial \theta_j}.
$$
However, in our case, $\Ktbt^{-1}$ and $\frac{\partial \mathbf K(\theta)}{\partial \theta_j}$ do not commute, and hence we need to find a way to compute $\frac{\partial \log(\Ktbt)}{\partial \theta_j}$ efficiently. For this, we use a well known result from the theory of matrix functions \cite{higham2008functions,al2011computing}: The directional derivative $L_{\mathrm{log}}(\Ab,\mathbf{E})$ of the matrix logarithm $\log(\Ab)$ in the direction $\mathbf{E}$ can be computed via 
$$
\log\left(
\begin{bmatrix}
\Ab&\mathbf{E}\\
0&\Ab
\end{bmatrix}
\right)=
\begin{bmatrix}
\log(\Ab)&L_{\mathrm{log}}(\Ab,\mathbf{E})\\
0&\log(\Ab)
\end{bmatrix}.
$$
Choosing the direction
$
E= \frac{\partial \mathbf{K}(\theta)}{\partial \theta_j},
$
we obtain
$$
\frac{\partial \log(\Ktbt)}{\partial \theta_j} = L_{\mathrm{log}}\left(\Ktbt,\frac{\partial \mathbf{K}(\theta)}{\partial \theta_j}\right).
$$
This can can be computed efficiently using the nonsymmetric TT-Krylov method and using the fact that
\begin{equation}
\label{approxFAB}
\begin{bmatrix}
\textbf{z}_i^{\mathsf T} \ 0
\end{bmatrix}
\log\left(
\begin{bmatrix}
\Ktbt&\frac{\partial \mathbf{K}(\theta)}{\partial \theta_j}\\
0&\Ktbt
\end{bmatrix}
\right)
\begin{bmatrix}
0\\
\textbf{z}_i
\end{bmatrix}
=\textbf{z}_i^{\mathsf T}L_{\mathrm{log}}\left(\Ktbt,\frac{\partial \mathbf{K}(\theta)}{\partial \theta_j}\right)\textbf{z}_i.
\end{equation}
We note that $\begin{bmatrix}
\textbf{z}_i^{\mathsf T} \ 0
\end{bmatrix}$ and $\begin{bmatrix}
0\\
\textbf{z}_i
\end{bmatrix}$ are TT tensors of order $D+1$ and with TT rank 1. Altogether,
we obtain an efficient way to compute the derivative of the trace estimator used in~\eqref{eq::NLL2}.

\paragraph{Derivatives for the tensor parameters}
For the above computations, we need to evaluate the derivative $\frac{\partial \mathbf{K}(\theta)}{\partial \theta_j}$ in TT format. We recall that
\begin{align*}
\mathbf{K} = \sum_{r=1}^R \bigotimes_{d=1}^D \mathbf{K}^{(d)}_r, \text{ where }  \mathbf{K}^{(d)}_r:= k^{(d)}_r(\textbf{X}^{(d)},\textbf{X}^{(d)}) \in \mathbb{R}^{n_d \times n_d},
\end{align*}
with the $1$-dimensional covariance matrices given by squared-exponential covariance functions and the signal variance shifted to the first covariance matrix, see \eqref{eq::single_covs}.

Hence, we can compute
\begin{align}\label{FullMat_deriv}
\frac{\partial\mathbf{K}}{\partial\ell_j^{(i)}} &= \frac{\partial\sum_{r=1}^R \bigotimes_{d=1}^D \mathbf{K}^{(d)}_r}{\partial\ell_j^{(i)}} =\sum_{r=1}^R \frac{\partial\bigotimes_{d=1}^D \mathbf{K}^{(d)}_r}{\partial\ell_j^{(i)}}\notag\\
&=\sum_{r=1}^R \sum_{d=1}^D \frac{\partial \mathbf{K}^{(d)}_r}{\partial\ell_j^{(i)}} \otimes \left(\bigotimes_{k \neq d} \mathbf{K}^{(k)}_r \right) \text{ (property of Kronecker product) }\notag\\
&= \frac{\partial \mathbf{K}_j^{(i)}}{\partial\ell_j^{(i)}} \otimes \left(\bigotimes_{d \neq i} \mathbf{K}^{(d)}_j \right),\text{ as } \frac{\partial \mathbf{K}^{(d)}_r}{\partial\ell_j^{(i)}} = 0 \text{ if } d \neq i \text{ or } r \neq j.
\end{align}
And similarly, for $r \in \{1,\ldots,R\}$:
\begin{align}\label{FullMat_deriv2}
\frac{\partial\mathbf{K}}{\partial\sigma_{f(r)}} = \frac{\partial \mathbf{K}^{(1)}_r}{\partial\sigma_{f(r)}} \otimes \left(\bigotimes_{d \neq 1} \mathbf{K}^{(d)}_r \right).
\end{align}
Here again, $\bigotimes_{d \neq i} \mathbf{K}^{(d)}_r$ can efficiently be computed as the Kronecker product of TT tensors. The matrix derivative $\frac{\partial \mathbf{K}_r^{(d)}}{\partial\theta}$ is now easy to compute for $\theta$ any length scale or signal variance parameter. As all covariance matrices result from squared exponential kernels, we use the log-transformed derivatives given by \eqref{eq::RBF_derivs} pointwise and obtain
\begin{align}\label{sigma_deriv}
\frac{\partial k_r^{(1)}(x,x^\prime)}{\partial\log(\sigma_{f(r)})} =2 k_r^{(1)}(x,x^\prime),
\end{align}
and
\begin{align}\label{ell_deriv}
\frac{\partial k_r^{(d)}(x,x^\prime)}{\partial\log(\ell_r^{(d)})} =\frac{(x-x^\prime)^2}{(\ell_r^{(d)})^2}k_r^{(d)}(x,x^\prime).
\end{align}
This computation has to be performed for each pair of inputs $x,x^\prime \in \textbf{X}^{(d)}$ to construct $\frac{\partial \mathbf{K}_r^{(1)}}{\partial\sigma_{f(r)}}$ and $\frac{\partial \mathbf{K}_r^{(d)}}{\partial\ell_r^{(d)}}$ for all possible $d \in \{1,\ldots,D\}$ and $r \in \{1,\ldots,R\}$. When the matrix of pointwise squared differences is also given in TT format, $\frac{\partial \mathbf{K}_r^{(d)}}{\partial\ell_r^{(d)}}$ can be computed from \eqref{ell_deriv} by a Hadamard product of TT matrices. 

We summarize the computation of the gradient in Alg.~\ref{alg:grad} and we now have all the ingredients to call {\matlab}'s $\mathtt{fminunc}$.

\begin{algorithm}[H]
	\caption{Algorithm for the computation of the gradient.}\label{alg:grad}
	\begin{algorithmic}[1]
		\Require $\theta = \{\sigma_{f(r)},\ell^{(d)}_r\}_{r = 1,\ldots,R,\, d = 1,\ldots,D}, p, \mathrm{kryltol}, \mathrm{amentol}$
		\Ensure $\nabla \hat f(\theta)$
		\For {$j = 1:\mathtt{length(\theta)}$}
		\State {$t = 0$}
		\For {$i=1:p$}
		\State {$\textbf{z}_i = \textbf{z}_i^{(1)} \otimes \cdots \otimes \textbf{z}_i^{(D)}, \; \textbf{z}_i^{(d)} \sim \mathcal U([-1,1]^{n_d}), \; d = 1,\ldots,D$}
		\State {$t = t + \begin{bmatrix}
\textbf{z}_i^{\mathsf T} \ 0
\end{bmatrix} \mathtt{tt\_krylov}\left(\mathtt{logm},\begin{bmatrix}
\Ktbt&\frac{\partial \mathbf{K}(\theta)}{\partial \theta_j}\\
0&\Ktbt
\end{bmatrix},
\begin{bmatrix}
0\\
\textbf{z}_i
\end{bmatrix},\mathrm{kryltol}\right)$}
		\EndFor	
		\State {$\mathbf \alpha = \mathtt{amen(\Ktbt,\mathbf y_{train},\mathrm{amentol})}$}
		\State {$\nabla \hat f(\theta)_j = t/p + \mathbf \alpha^{\mathrm T} \frac{\partial \mathbf{K}(\theta)}{\partial \theta_j} \mathbf \alpha$}
		\EndFor
	\end{algorithmic}
\end{algorithm}

\section{Numerical experiments}\label{sec:NumExp}

We test our method on two synthetic cases: In one, we generate a random trigonometric function in $3$ dimensions. In the second example, we draw a random sample from a Gaussian process with given covariance kernel. In both cases, we learn the parameters of our Gaussian process on some training points and evaluate the results on a range of test data. This is done for an increasing number of probe vectors in the trace estimation.

We start the construction of the test problem by generating a tensor ${\mathbf X_{train} \in \mathbb R^{N \times N \times N}}$ as a 3-dimensional grid in $[-1,1]^3$:
\begin{equation}
\mathbf X_{train} = \mathbf x_1 \otimes \mathbf x_2 \otimes \mathbf x_3,
\end{equation}
where the points are equally spaced via the construction
\begin{equation}
\mathbf x_1 = \mathbf x_2 = \mathbf x_3 = [-1:\frac{2}{N-1}:1].
\end{equation}
This tensor contains the training points. As test points, we generate the complimentary grid $\mathbf X_{test} \in \mathbb R^{N-1 \times N-1 \times N-1}$ 
\begin{equation}
\mathbf X_{test} = \mathbf {\tilde x_1} \otimes \mathbf {\tilde x_2} \otimes \mathbf {\tilde x_3},
\end{equation}
with
\begin{equation}
\mathbf {\tilde x_1} = \mathbf {\tilde x_2} = \mathbf {\tilde x_3} = [-1 + \frac{1}{N-1} :\frac{2}{N-1} : 1 - \frac{1}{N-1}].
\end{equation}
The right hand side is then given as $\mathbf y_{train} \in \mathbb R^{N \times N \times N}$ and $\mathbf y_{test} \in \mathbb R^{N-1 \times N-1 \times N-1}$, respectively.

\subsection{Trigonometric function with random coefficients}

In the first experiment, we set $N = 21$ and generate a random tensor $\mathbf R \in \mathbb R^{3 \times 3 \times 2}$ with entries uniformly distributed in the interval $[0,1]$. The training labels are then given as
\begin{equation}
\begin{aligned}
\mathbf y_{train}(i,j,k) = \sum_{r=1}^3 \sin(&\pi \mathbf R(r,1,1) \mathbf x_1(i) + \frac{\pi}{2} \mathbf R(r,1,2)) \sin(\pi\mathbf  R(r,2,1) \mathbf x_2(i) \\
+ &\frac{\pi}{2} \mathbf R(r,2,2)) \sin(\pi \mathbf R(r,3,1) \mathbf x_3(i) + \frac{\pi}{2} \mathbf R(r,3,2)),
\end{aligned}
\end{equation}
and similarly for the test points $\mathbf y_{test}$ with the same parameter tensor $\mathbf R$. We add the random noise with $\sigma = 0.01 $ and the resulting tensors are then converted into TT format with rounding error $10^{-8}$.

\begin{figure}[!b]
\centering
\includegraphics[width=.5\textwidth]{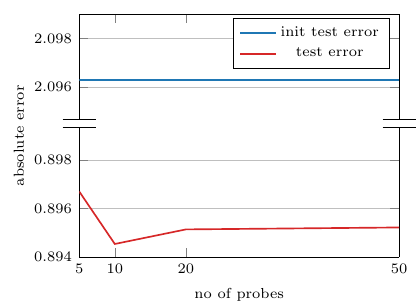}
\caption{ Comparison of the test error of the Gaussian process with the initial hyperparameters (blue) and with optimized hyperparameters (red) plotted against an increasing number of probe vectors for the trace estimation. }
\label{fig:funcerr}
\end{figure}

Since the training points are generated by a sum of 3 products of univariate trigonometric functions, we will attempt to reconstruct this function using a Gaussian process with a Kronecker-sum kernel of rank $R = 3$ and order $D = 3$. 
We initialize our method with parameters $\sigma_{f(r)} = 1 + \varepsilon, r = 1,2,3$ and $\ell_r^{(d)} = 0.1 + \varepsilon, r = 1,2,3, d = 1,2,3$ where $\varepsilon$ is random noise in $\mathcal N(0,0.005)$ and different for each parameter.

We set the tolerance for the TT-Krylov method and the \textsc{AMEn} subsolver to $10^{-6}$ and we run our method using the LBFGS solver in {\matlab}'s \texttt{fminunc} routine (the tolerance for the norm of the gradient was set to $10^{-10}$).
We test different numbers of probe vectors in the trace estimation and report on the results.

In Fig.~\ref{fig:funcerr}, we show the absolute Frobenius error of the posterior mean tensor $\bm{\mu}_*$ to the test tensor $\mathbf y_{test}$. We can see that in all cases, the error is significantly reduced as compared to the random initialization. However, the error does not visibly improve for different numbers of probe vectors $p$.

\begin{figure}[!t]
\subfloat[Trace estimation error\label{subfig:trace}]{%
\includegraphics[width=0.5\textwidth]{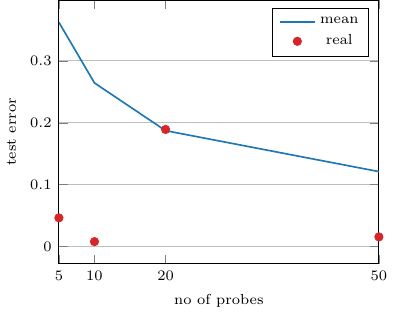}
}
\hfill
\subfloat[Gradient approximation\label{subfig:grad}]{%
\includegraphics[width=0.5\textwidth]{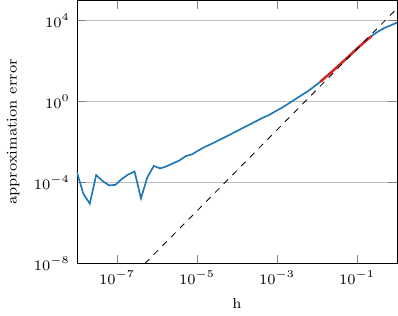}
}
\caption{ (a) Error of the trace estimation based on the  number of probes with the mean approximation (blue line) and the actual value (red). (b) Approximation error of the cost function $\hat f$ by its first order approximation using the computed gradient, for different step sizes $h$. The dashed line represents the expected approximation error of order 1, the red part highlights the interval of step sizes where this theoretical error is matched.}
\label{fig:functrace}
\end{figure}

\begin{figure}[!b]
\centering
\includegraphics[width=\textwidth]{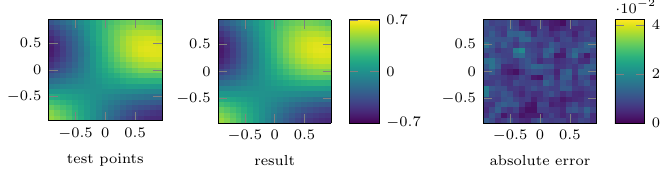}
\caption{Slice of a data tensor on the left for the test data and in the middle the same slice of the posterior mean using the optimized parameters. The error is shown on the right.}
\label{fig:funcheat}
\end{figure}

This is most likely due to the random nature of the trace estimation. In Fig.~\ref{fig:functrace}, we show the relative error of the trace estimation in the first iteration. We can see that this error varies significantly from the mean. This explains the fact that we can get good approximation errors with only a small number of probes. On the right hand side, we see the approximation error of the cost function compared to its first-order Taylor approximation using our computed gradient. One expects that this approximation error grows linearly with the step size $h$. We can see that this is indeed the case for relevant step sizes $h \in [0.01,0.5]$ but the approximation is less exact for smaller or larger step sizes due to the inaccurate computation of the gradient. 

Finally, in Fig.~\ref{fig:funcheat}, we show a comparison of one slice of the test tensor and the posterior mean for 10 probes (as this gave the best test error). We chose the $10$th slice and compare $\mathbf y_{test}(:,10,:)$ with $\bm{\mu}_*(:,10,:)$ as well as the absolute error of the two. We can see that in this experiment, the interpolation works very well and the point-wise errors are small.

\subsection{Random sample from Gaussian process}

In the second example, we repeat the previous experiment but with a different right hand side $\mathbf y_{train}$ and also different test points $\mathbf y_{test}$. We take the same grid $\mathbf X$ but generate these tensors by drawing them randomly from the Gaussian process with Kronecker covariance kernel $\mathbf{K} = \sum_{r=1}^3 \bigotimes_{d=1}^3 \mathbf{K}^{(d)}_r$, using the parameters
\begin{align*}
\sigma_{f(1)} &= 1, &\sigma_{f(2)} &= 0.1, & \sigma_{f(3)} &= 0.01, \\
\ell_1^{(1)} &= 0.06, &\ell_2^{(1)} &= 0.2,  & \ell_3^{(1)} &= 0.3, \\
\ell_1^{(2)} &= 0.05, &\ell_2^{(2)} &= 0.19, & \ell_3^{(2)} &= 0.4, \\
\ell_1^{(3)} &= 0.04, &\ell_2^{(3)} &= 0.21, & \ell_3^{(3)} &= 0.5.
\end{align*}
The idea is here that short-range interactions ($\ell_r^{(d)} < 1/N$) outweigh the longer interactions in the Kronecker terms 2 and 3.

The training tensor is produced by drawing two random tensors ${\mathbf U,\mathbf V \in \mathbb R^{(2N-1) \times (2N-1) \times (2N-1)}}$ with entries drawn from the standard normal distribution, converting them to the TT format with rounding error $10^{-8}$ and computing
\begin{equation}
\mathbf y_{full} = \mathbf K^{1/2} \mathbf U + \sigma \mathbf V,
\end{equation}
where $\sigma = 0.01$ is the added observational noise. The matrix square root is computed using our Krylov method for TT tensors and $g(\cdot) = \sqrt{\cdot}$. The resulting tensor is then rounded again with accuracy $10^{-8}$. We then compute the training and test tensors by
\begin{align}
\mathbf y_{train}(i,j,k) &= \mathbf y_{full}(2i-1,2j-1,2k-1), \qquad &i,j,k = 1,\ldots,N, \\
\mathbf y_{test}(i,j,k) &= \mathbf y_{full}(2i,2j,2k), &i,j,k = 1,\ldots,N-1.
\end{align}

We again reconstruct this function using a Gaussian process with covariance kernel of rank $R = 3$ and order $D = 3$ and we initialize the method with the same values used above. All other parameters remain the same as well.

\begin{figure}[!t]
\centering
\includegraphics[width=.5\textwidth]{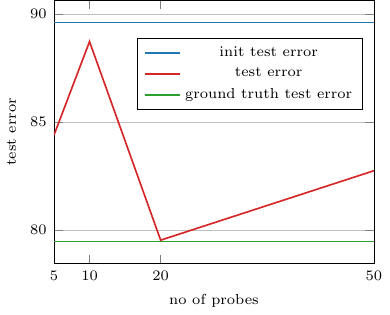}
\caption{Comparison of the test error of the Gaussian process with the initial hyperparameters (blue), optimized hyperparameters (red), and with the ground truth parameters (green), plotted against an increasing number of probe vectors for the trace estimation.}
\label{fig:randerr}
\end{figure}

In Fig.~\ref{fig:randerr}, we again show the absolute Frobenius error of the posterior mean tensor $\bm{\mu}_*$ to the test tensor $\mathbf y_{test}$. Apart from the error for the posterior mean using the initial parameters we also show the test error for the (known) ground truth parameters set above. We can see that the test error of our method fluctuates. In the case of 20 probes, we get about the same test error as for the ground truth.

However, this experiment is clearly harder. In Fig.~\ref{fig:randheat}, we show the 10th slice of the computed tensor in the case of 20 probes, as compared to the test tensor. In this case, we see that the error is much larger, and the resulting tensor is only a blurry approximation of the original function. We cannot expect a much better result, since the test error is the same also for the ground truth. This is due to the random nature of the problem and it likely being more difficult to approximate with a low-rank method given the little smoothness that the data possesses.

\begin{figure}[!b]
\centering
\includegraphics[width=\textwidth]{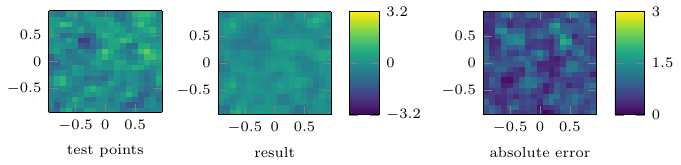}
\caption{Slice of the tensor for test points (right), optimized setup (middle) and error (right).}
\label{fig:randheat}
\end{figure}

\section{Conclusions}
In this paper we proposed a numerical scheme for the optimization of the hyperparameters of a Gaussian process model for tensor-valued data that are approximated in the tensor train format. The method relies on the solution of an equation on the \textsc{AMEn} solver and we derived a matrix function analogue for the evaluation of both the objective function and its derivative. The estimation of the trace was done using a matrix logarithm and for its approximation Krylov methods and rational Krylov methods are proposed. For the derivative we require the evaluation of a block matrix function to get the Frechet derivative. 
We illustrate that the method indeed approximates the theoretical gradients well for the relevant range of parameters and show for two synthetic examples that we produce good approximations.

\bmhead{Acknowledgments}
The authors would like to thank David Bindel for the insightful discussions. M.S. acknowledges discussions with Kim Batselier on the use of tensor methods for Gaussian processes. The research of J.K.~was partially funded by the Deutsche Forschungsgemeinschaft (DFG, German Research Foundation)---Project-ID 318763901 - SFB1294. M.P.~was partially funded by the DFG – Projektnummer 448293816.

\bibliography{refs}

\end{document}